\documentstyle[12pt]{article}
\if@twoside  m
\else
\fi
   \newcommand{\ie}{{\em i.e.}}
   \newcommand{\eg}{{\em e.g.}}
   \newcommand{\rhs}{{\em rhs }}
   \newcommand{\lhs}{{\em lhs }}
   \newcommand{\LL}{{\cal L}}
   \newcommand{\OO}{{\cal O}}
   \newcommand{\R}{I\!\!R}
   \newcommand{\supp}{{\rm supp\,}}
   \newcommand{\re}{{\rm Re\,}}
   \newtheorem{claim}{Claim}
   \newtheorem{theorem}[claim]{Theorem}
\begin{document}
\vspace*{20mm}

\noindent
{\Large\bf Bound states in a locally deformed waveguide: \\ the
critical case}
\vspace{5mm}

\noindent
{\bf P.~Exner} \\
{\small Nuclear Physics Institute, Academy of Sciences,
25068 \v{R}e\v{z} near Prague, \\ and Doppler Institute,
Czech Technical University, B\v rehov\'a 7, 11519 Prague, \\
Czech Republic, \\ \em exner@ujf.cas.cz} \\
and \\
{\bf S.A.~Vugalter} \\{\small Nuclear Physics Institute,
Academy of Sciences, 25068 \v{R}e\v{z} near Prague,
Czech Republic, \\ and Radiophysical Research Institute, B.
Pecherskaya 25/14, 603600 Nizhni Novgorod, Russia}
\vspace{10mm}

\noindent
{\small {\bf Abstract.} We consider the Dirichlet Laplacian for a
strip in $\,\R^2$ with one straight boundary and a width
$\,a(1+\lambda f(x))\,$, where $\,f\,$ is a smooth function of a
compact support with a length $\,2b\,$. We show that in the critical
case, $\,\int_{-b}^b f(x)\,dx=0\,$, the operator has no bound states
for small $\,|\lambda|\,$ if $\,b<(\sqrt{3}/4)a\,$.  On the other
hand, a weakly bound state exists provided $\,\|f'\|< 1.56
a^{-1}\|f\|\,$; in that case there are positive $\,c_1, c_2\,$ such
that the corresponding eigenvalue satisfies $\,-c_1\lambda^4\le
\epsilon(\lambda)- (\pi/a)^2 \le -c_2\lambda^4\,$ for all
$\,|\lambda|\,$ sufficiently small. }
\vspace{10mm}

\section{Introduction}

Quantum waveguides, \ie, systems in which a quantum mechanical
particle is confined to a tubular region, attracted a wave of
attention induced by the progress of the ``mesoscopic" physics ---
some reference and a guide to further reading can be found in
\cite{DE,DES}. Apart of a natural physical interest, some interesting
mathematical problems arose in this connection.

One of them concerns the existence of bound states which appear if a
Dirichlet tube of a constant cross section is locally deformed, \eg,
bent --- see \cite{ES,DE} and references therein.  Another mechanism
producing bound states in straight tubes is based on tube protrusions
\cite{ISY,AS}. It has been treated rigorously in a recent paper by
Bulla et al. \cite{BGRS} where the asymptotic behavior of the
eigenvalue for a gentle deformation was found.

In distinction to bent tubes a local variation of the cross section
can yield both an attractive and repulsive effective interaction;
it is easy to see that constrictions produce no bound states. A more
subtle question is what happens if the profile modification along the
tube includes both an expansion and squeezing. It is demonstrated
in \cite{BGRS} that the ``averaged" cross--section variation matters:
a bound state exists in a slightly deformed tube provided the added
volume is positive while in the opposite case it is absent. In the
former case, the asymptotic behavior of the gap is governed by the
square of the deformation parameter $\,\lambda\,$.

In this letter we address the critical case left out in the mentioned
paper, namely the situation when the volume change is zero. The
authors pointed out that bound states might exist in view of the 
analogy with one--dimensional Schr\"odinger operators
\cite{BGS,Si,Kl}, however, they remarked at the same time that the
effect of the additional second--order terms in the Hamiltonian is
not apriori obvious.  For the sake of simplicity we restrict
ourselves to the basic setup of \cite{BGRS}: we consider a planar
strip obtained be deforming locally {\em one} boundary of
$\,\Omega_0:= \R\times [0,a]\,$. We show that the answer depends on
the shape of the deformation. If the latter is localized into an
interval sufficiently small with respect to tube width $\,a\,$, there
is no bound state. On the other hand, a bound state exists if the
deformation is smeared enough, and the leading behavior of the gap is
given in this case by the {\em fourth} power of the parameter
$\,\lambda\,$. 

\section{The results}

We consider the Dirichlet Laplacian $\,-\Delta^D_{\Omega_\lambda}\,$
on a deformed strip
   \begin{equation} \label{strip}
\Omega_\lambda\,:=\, \{\, (x,y)\in\R^2:\; 0<y<a(1+\lambda
f(x))\,\}\,,
   \end{equation}
where $\,f\in C_0^{\infty}(\R)\,$ is a given function; in distinction
to \cite{BGRS} we do not put $\,a=1\,$. Our main results are the
following:

   \begin{theorem} \label{nonexistence}
Suppose that $\,f\in C_0^{\infty}(\R)\,$ is such that $\,\supp
f\subset [-b,b]\,$ and
   \begin{equation} \label{zero mean}
\int_{-b}^b f(x)\,dx\,=\,0\,;
   \end{equation}
then the discrete spectrum of $\,-\Delta^D_{\Omega_\lambda}\,$ is
empty for all sufficiently small $\,|\lambda|\,$ provided
   \begin{equation} \label{nonexistence condition}
a\,>\, {4\over \sqrt 3}\:b \,.
   \end{equation}   \end{theorem}

   \begin{theorem} \label{existence}
Under the same assumptions, $\,-\Delta^D_{\Omega_\lambda}\,$ has an
isolated eigenvalue $\,\epsilon(\lambda)\,$ for any nonzero
$\,\lambda\,$ with $\,|\lambda|\,$ small enough provided 
   \begin{equation} \label{existence condition}
{\|f'\|^2\over \|f\|^2}\,<\, \left( \pi\over a\right)^2\, {6\over
9+\sqrt{90+12\pi^2}}\,.
   \end{equation}
In that case there are positive $\,c_1,\,c_2\,$ such that
   \begin{equation} \label{asymptotics}
-c_1\lambda^4 \,\le\, \epsilon(\lambda)-\, \left(\pi\over
a\right)^2 \,\le\, -c_2\lambda^4\,.
   \end{equation}   \end{theorem}

\section{Nonexistence of bound states}

By \cite{BGRS}, the operator under consideration is unitarily
equivalent to
   \begin{equation} \label{Hamiltonian}
H_{\lambda} =\, -\Delta^D_{\Omega_0} +\lambda A_1 -\lambda^2 A_2
   \end{equation}
on $\,L^2(\Omega_0)\,$ corresponding to the straightened strip
$\,\Omega_0:= \R\times [0,a]\,$, where
$$
A_1\,=\, \sum_{j=1}^5 A_{1j}\,, \qquad A_2\,=\, \sum_{j=1}^6 A_{2j}
$$
with
   \begin{eqnarray} \label{first order}
&& A_{11}= 2f(x)\partial^2_y\,, \quad A_{12}= yf''(x)\partial_y\,,
\quad A_{13}= 2yf'(x) \partial^2_{xy}\,, \nonumber \\
&& A_{14}= f'(x)\partial_x\,, \quad A_{15}= {1\over 2}\,f''(x)\,,
   \end{eqnarray}
and
   \begin{eqnarray} \label{second order}
&& \hspace{-20pt} A_{21}= {3f(x)^2+2\lambda f(x)^3+y^2f'(x)^2 \over
(1+\lambda f(x))^2}\, \partial^2_y\,, \quad
A_{22}= \left( {yf(x)f''(x)\over 1+\lambda f(x)} + {3yf'(x)^2\over
(1+\lambda f(x))^2} \right) \partial_y\,, \nonumber \\ \nonumber \\
&& \hspace{-20pt} A_{23}= {2yf(x)f'(x)\over 1+\lambda f(x)}\,
\partial^2_{xy}\,, \qquad A_{24}= {f(x)f'(x)\over 1+\lambda f(x)}\,
\partial_x\,, \\ \nonumber \\
&& \hspace{-20pt} A_{25}= {f(x)f'(x)\over 2(1+\lambda f(x))}\,,
\qquad A_{26}= {3f'(x)^2\over 4(1+\lambda f(x))^2}\,. \nonumber
   \end{eqnarray}
Since $\,\inf \sigma_{ess}(H_\lambda)= (\pi/a)^2$,
Theorem~\ref{nonexistence} will be proven if we demonstrate that
$\,(H_\lambda\psi, \psi)\ge (\pi/a)^2$ holds for all $\,\psi\,$ from
a suitable dense set, say, $\,C_0^2(\Omega_0)\,$. Such a function can
be always written as
   \begin{equation} \label{decomposition}
\psi(x,y)\,=\, G(x,y)+R(x,y)\,,
   \end{equation}
where $\,G(x,y)=\phi(x)\chi_1(y)\,$ and $\,R(x,\cdot)\perp \chi_1\,$
for all $\,x\in\R\,$; we use the symbol $\,\chi_n\,$ for the
normalized $\,n$--th transverse--mode eigenfunction,
   \begin{equation} \label{transverse mode}
\chi_n(y)\,:=\, \sqrt{2\over a}\, \sin\left(\pi ny\over a\right)\,.
   \end{equation}
The smooth function $\,\phi\,$ can be written as $\,\phi(x)= \alpha+
g(x)\,$, where $\,\alpha:= \phi(-b)\,$. In view of 
(\ref{Hamiltonian}), the quadratic form to be estimated can be
expressed as
   \begin{eqnarray} \label{quadratic form}
(H_\lambda\psi,\psi) \! &=& \! -(\Delta^D_{\Omega_0}G,G)
-(\Delta^D_{\Omega_0}R,R) \nonumber \\
&& +\lambda(A_1G,G) +\lambda(A_1R,R) +2\lambda\,\re(A_1G,R) \\
&& -\lambda^2(A_2G,G) -\lambda^2(A_2R,R) -2\lambda^2\re(A_2G,R)\,.
\nonumber
   \end{eqnarray}
Let us begin with the terms linear in $\,\lambda\,$. With an abuse of
notation, $\,A_1(1\times\chi_1)= A_1\chi_1\,$, we write
$$
(A_1G,G)\,=\, \alpha^2(A_1\chi_1,\chi_1) +\alpha(A_1\chi_1,g\chi_1)
+(A_1g\chi_1,g\chi_1)\,.
$$
The first term at the \rhs is zero in view of (\ref{zero mean}) and
$\,\int_{-b}^b f''(x)\,dx= f'(b)\!-\!f'(-b)=0\,$. The second one is
$$
-2\alpha\,\left(\pi\over a \right)^2 (f,g)_{L^2(-b,b)} -{\alpha\over
2}\, (f',g')_{L^2(-b,b)} -\alpha (y\chi'_1,\chi_1)_{L^2(0,a)}
(f',g')_{L^2(-b,b)}\,,
$$
where we have used integration by parts together with $\,-\chi''_1=
(\pi/a)^2 \chi_1\,$; evaluating the inner product we see that the
last two terms cancel. Finally, using the fact that $\,f,f',f''$ are
bounded on $\,[-b,b]\,$, we arrive at the bound
$$
|(A_1g\chi_1,g\chi_1)|\,\le\, C\left(\|g\|^2_{L^2(-b,b)}
+\|g'\|^2_{L^2(-b,b)} \right)\;;
$$
here and in the following $\,C\,$ is an unspecified positive constant
which assume different values in different expressions. In fact, the
derivative norm alone may be used here, because
$\,g(-b)=0\,$ implies easily
   \begin{equation} \label{embedding}
\|g\|^2_{L^2(-b,b)} \,\le\, (2b)^2 \|g'\|^2_{L^2(-b,b)}\,.
   \end{equation}
Together we have
$$
\lambda(A_1G,G)\,\ge\, -2\alpha\,\left(\pi\over a \right)^2
(f,g)_{L^2(-b,b)} -\lambda C\|g'\|^2_{L^2(-b,b)}\,.
$$
In a similar way the boundedness of $\,f,f',f''$ together with
integration by parts and the Schwarz inequality yield
$$
|(A_1R,R)|\,\le\, C\|R\|^2_{W_2^1(\Omega_b)}\,,
$$
where $\,\Omega_b:= [-b,b]\times[0,a]\,$. The same argument applies
to the ``$\,\alpha$--indepedent part" of the mixed term,
$$
|(A_1g\chi_1,R)|\,\le\, C \left(\|g'\|^2_{L^2(-b,b)}+
\|R\|^2_{L^2(\Omega_b)} \right)\,.
$$
The remaining term needs more attention. Since $\,R\,$ is smooth by
assumption, it may be represented pointwise as
   \begin{equation} \label{decomposition 2}
R(x,y)\,=\, \sum_{n=2}^\infty r_n(x)\chi_n(y)\,,
   \end{equation}
so
   \begin{equation} \label{mixed term 1}
\alpha(A_1\chi_1,R)\,=\, -\alpha\, \sum_{n=2}^\infty
(f',g')_{L^2(-b,b)} (y\chi'_1,\chi_n)_{L^2(0,a)}\,.
   \end{equation}
The last inner product equals $\,(-1)^n 2n/(n^2\!-\!1)\,$; we introduce
   \begin{equation} \label{K}
K\,:=\, \left(\sum_{n=2}^\infty \left( 2n\over n^2\!\!-\!1 \right)^2
\right)^{1/2}\,.
   \end{equation}
The last sum can be evaluated \cite[Sec.~5.1]{PBM} in terms of the
di-- and trigamma functions $\,\psi\,$ and $\,\psi'$, respectively, as
$$
4\,\sum_{n=1}^\infty {(n+1)^2\over n^2(n+2)^2}\,=\, \pi^2-\,
{31\over 4}\,+\, 2\big(\psi(3)+\gamma\big)\,-\,4\psi'(3)\,,
$$
and since $\,2(\psi(3)\!+\!\gamma)=3\,$ and $\,4\psi'(3)=
4\zeta(2)\!-\!5\,$ by \cite[Chap.~6]{AbS}, where $\,\zeta\,$ is the
Riemann zeta function and $\,\gamma\,$ is the Euler's constant, we
obtain  
   \begin{equation} \label{K explicit}
K\,=\, \sqrt{\,{\pi^2\over 3}+{1\over 4}}\,.
   \end{equation}
The term (\ref{mixed term 1}) in question may be then estimated in
modulus by
   \begin{eqnarray*}
&& \alpha\left(|f'|, \sum_{n=2}^\infty |r'_n| {2n\over n^2\!-\!1}
\right)_{L^2(-b,b)} \,\le\, \alpha K\left(|f'|, \left(
\sum_{n=2}^\infty |r'_n|^2 \right)^{1/2}\right)_{L^2(-b,b)} \\ \\
&& \le\,{1\over 2}\, \alpha^2 K^2 \lambda (1\!+\!\tilde c\lambda)
\|f'\|^2_{L^2(-b,b)} +\, {1\over 2\lambda(1\!+\!\tilde c\lambda)}
\|R_x\|^2_{L^2(\Omega_b)} \,,
   \end{eqnarray*}
where $\,R_x:= \partial R/\partial x\,$ and $\,\tilde c\,$ is a
positive number to be specified later; in the last step we have used
the Schwarz inequality.

Next, we pass to the quadratic terms. By the boundedness of
$\,f,f',f''\,$ we have
$$
\lambda^2|(A_2R,R)|\,\le\, C\lambda^2 \|R\|^2_{W_2^1(\Omega_b)}\;;
$$
the same property in combination with the Schwarz inequality and the
estimate (\ref{embedding}) gives for the two parts of the mixed term
   \begin{eqnarray*}
\lambda^2\alpha |(A_2\chi_1,R)| \! &\le & \! C\alpha^2\lambda^3+\,
{\lambda\over 2}\, \|R\|^2_{L^2(\Omega_b)} \\ \\
\lambda^2 |(A_2 g\chi_1,R)| \! &\le & \! C \lambda^3
\|g'\|^2_{L^2(-b,b)} + \, {\lambda\over 2}\, \|R\|^2_{L^2(\Omega_b)}\,.
   \end{eqnarray*}
The remaining quadratic term which refers to the lowest
transverse--mode component equals
$$
\lambda^2\left( \alpha^2 (A_2\chi_1,\chi_1)+ 2\alpha\,\re
(A_2g\chi_1,\chi_1)+ (A_2g\chi_1,g\chi_1) \right)\,.
$$
The last two terms at the \rhs are estimated as above
   \begin{eqnarray*}
\lambda^2\alpha |(A_2g\chi_1,g\chi_1)| \! &\le & \! C\lambda^2
\|g'\|^2_{L^2(-b,b)} \\ \\
\lambda^2 |(A_2 \chi_1,g\chi_1)| \! &\le & \! C \alpha^2\lambda^3 +
\, {\lambda\over 2}\, \|g\|^2_{L^2(-b,b)}\,,
   \end{eqnarray*}
while the first one can be evaluated from
   \begin{eqnarray*}
\lefteqn{(A_2\chi_1,\chi_1) \,=\, -3\left(\pi\over a\right)^2
\|f\|^2_{L^2(-b,b)} -\left(\pi\over a\right)^2\|f'\|^2_{L^2(-b,b)}
(y^2\chi_1,\chi_1)_{L^2(0,a)}} \\ \\
&& +2\|f'\|^2_{L^2(-b,b)} (y\chi'_1,\chi_1)_{L^2(0,a)} +\,{1\over
4}\|f'\|^2_{L^2(-b,b)}  +\OO(\lambda)\,.
   \end{eqnarray*}
Computing the inner products and using (\ref{K explicit}), we arrive
thus at the expression
$$
\lambda^2 (A_2\chi_1,\chi_1)\,=\, -\lambda^2\alpha^2
\left(3\left(\pi\over a\right)^2 \|f\|^2_{L^2(-b,b)} +K^2
\|f'\|^2_{L^2(-b,b)} \right) +\OO(\lambda^3) \,.
$$
Putting the results obtained up to now together, we find for the \rhs
of (\ref{quadratic form}) the bound
   \begin{eqnarray} \label{form lower bound}
\lefteqn{(H_\lambda\psi,\psi) \,\ge\, \|\nabla G\|^2 + \|\nabla R\|^2
+3\lambda^2\alpha^2 \left(\pi\over a\right)^2\|f'\|^2_{L^2(-b,b)}
-2\lambda\alpha \left(\pi\over a\right)^2 (f,g)_{L^2(-b,b)}}
\nonumber \\ \\
&& -(1\!+\!\tilde c\lambda)^{-1} \|R_x\|^2_{L^2(\Omega_b)}
-C\left(\lambda \|g'\|^2_{L^2(-b,b)} +\lambda
\|R\|^2_{W_2^1(\Omega_b)} +\alpha^2\lambda^3 \right)\,,
\phantom{AAAAAAA} \nonumber
   \end{eqnarray}
which is valid for all $\,|\lambda|\,$ small enough. It remains to
estimate the ``kinetic" terms. Using once more the decomposition
(\ref{decomposition 2}), we derive
$$
\|\nabla R\|^2\,\ge\, \|R_x\|^2_{L^2(\Omega_b)} +\left(2\pi\over
a\right)^2 \|R\|^2_{L^2(\Omega_b)}\;;
$$
the inequality
$$
\|\nabla R\|^2 -\lambda C\|R\|^2_{W_2^1(\Omega_b)} -(1\!+\!\tilde
c\lambda)^{-1} \|R_x\|^2_{L^2(\Omega_b)} \,\ge\, \left(2\pi\over
a\right)^2 \|R\|^2_{L^2(\Omega_b)}
$$
is then satisfied if
$$
{\tilde c-C\over 1\!+\!\tilde c\lambda}\,>\,0 \qquad {\rm and} \qquad
3\left(\pi\over a\right)^2 \!-\lambda C\,>\,0\,,
$$
\ie, for $\,\tilde c>C\,$ and all $\,|\lambda|\,$ small enough. In a
similar way, we obtain
$$
\|\nabla G\|^2\,\ge\, \|g'\|^2_{L^2(-b,b)} +\left(\pi\over a\right)^2
\|G\|^2_{L^2(\Omega_b)}\,.
$$
Without loss of generality we may assume $\,f\ne 0\,$. We insert the
gradient--term estimates into (\ref{form lower bound}) and estimate
$\,\|g'\|^2_{L^2(-b,b)}\,$ from below by $\,\|g\|^2_{L^2(-b,b)}\,$.
This can be done by the inequality (\ref{embedding}), however, the
latter is unnecessarily rough. Minimization of
$\,\|g'\|^2_{L^2(-b,b)}\,$ for fixed $\,g(-b)=0\,,\; g(b)\,$, and
$\,\|g\|^2_{L^2(-b,b)}\,$ is an isoperimetric problem
\cite[Sec.23F]{Re}; the corresponding Euler's equation is easily seen
to be solved by a multiple of $\,\sin \kappa(x\!-\!b)\,$. Taking a
minimum over $\,\kappa\,$, we get  
   \begin{equation} \label{embedding 2}
\|g'\|^2_{L^2(-b,b)} \,\ge\, \left(\pi\over 4b\right)^2
\|g\|^2_{L^2(-b,b)}\,.
   \end{equation}
Using this bound, we arrive at
   \begin{eqnarray*}
\lefteqn{(H_\lambda\psi,\psi)-\, \left(\pi\over a\right)^2
\|\psi\|^2 \,\ge\, \pi^2\:{1-C\lambda\over 16b^2}\,
\|g\|^2_{L^2(-b,b)}} \\ \\ && -\,2\alpha\lambda
\left(\pi\over a\right)^2 (f,g)_{L^2(-b,b)} +\lambda^2
\|f\|^2_{L^2(-b,b)} \left(3\alpha^2 \left(\pi\over a\right)^2-\,
{C\alpha^2 \lambda\over \|f\|^2_{L^2(-b,b)}} \right)\;;
   \end{eqnarray*}
the quadratic form at the \rhs is positive as long as
   \begin{equation} \label{discriminant}
\left(\pi\over a\right)^4\,<\, \pi^2\,
\left(1-C\lambda\over 16b^2 \right) \left(3 \left(\pi\over
a\right)^2-\, {C \lambda\over \|f\|^2_{L^2(-b,b)}} \right)\,.
   \end{equation}
If the inequality is satisfied for $\,\lambda=0\,$ the same is true
for $\,|\lambda|\,$ small enough; this yields the condition
(\ref{nonexistence condition}) and finishes thus the proof of
Theorem~\ref{nonexistence}.

\section{Existence of bound states}

Let us start from a lower bound to possible eigenvalues. A necessary
condition for their existence is that the condition
(\ref{discriminant}) is violated. It follows easily from the above
estimates that
   \begin{equation} \label{bs lower bound}
(H_\lambda\psi,\psi)-\,\left(\pi\over a\right)^2 \|\psi\|^2 \,\ge\,
\|G_x\|^2_{L^2(\Omega^c_b)} -\lambda^2\alpha^2 d_0 \|f\|^2
+\OO(\lambda^3)\,, 
   \end{equation}
where $\,\Omega^c_b:=\Omega_0\setminus\Omega_b\,$ is the complement
containing the strip tails, and
$$
d_0\,:=\, \left(4\pi b\over a^2\right)^2 -3\left(\pi\over
a\right)^2\;; 
$$
if a bound state exists this quantity has to be positive.
Consider a function $\,\psi\,$ for which the \lhs of (\ref{bs
lower bound}) is negative. We employ the inequality $\,\|\psi\|^2 \ge
\|G\|^2_{L^2(\Omega^{c,-}_b)}\,$, where $\,\Omega^{c,-}_b\,$ is the
left tail, and an analogous bound for
$\,\|G_x\|^2_{L^2(\Omega^c_b)}\,$. The functional 
$$
\LL(G)\,:=\, {\|G_x\|^2_{L^2(\Omega^{c,-}_b)} -\lambda^2\alpha^2 d_0
\|f\|^2 \over \|G\|^2_{L^2(\Omega^{c,-}_b)}} \,=\,
{\|\phi'\|^2_{L^2(-\infty,-b)} -\lambda^2\alpha^2 d_0
\|f\|^2 \over \|\phi\|^2_{L^2(-\infty,-b)}}
$$
defined on functions with fixed $\,\phi(-b)=\alpha\,$ assumes its
extremum for $\,G_0(x,y)= \alpha\, e^{\kappa(x+b)} \chi_1(y)\,$
with $\,\kappa=\lambda^2 d_0\|f\|^2$; the minimum value is
$\,-\lambda^4 d_0^2\|f\|^4$. Consequently, we get
$$
{(H_\lambda\psi,\psi)\over \|\psi\|^2} -\,\left(\pi\over a\right)^2
\,\ge\, -d_0^2 \|f\|^4 \lambda^4 +\OO(\lambda^5)\,,
$$
\ie, the lower bound (\ref{asymptotics}) of Theorem~\ref{existence}.

The rest of the argument is easier. To check the existence claim
of Theorem~\ref{existence}, it is sufficient to use a suitable trial
function. We choose
   \begin{equation} \label{trial function}
\psi_\lambda(x,y)\,=\, \big(1\!+\!\lambda\eta f(x)\big)
\chi_1(y)\,,
   \end{equation}
where $\,\eta\,$ is a parameter to be determined; this allows us
to employ the above estimates with $\,R=0\,\; \alpha=1\,$, and
$\,g=\lambda\eta f\,$. In particular, we have
   \begin{eqnarray*}
(H_\lambda\psi,\psi)\,-\,\left(\pi\over a\right)^2
\|\psi\|^2  \!&\le&\! \|\psi_x\|_{L^2(\Omega_b^c)} \\ \\
\!&+&\! \lambda^2 \left\{(3\!-\!2\eta) \left(\pi\over a\right)^2
\|f\|^2 + (\eta^2\!+K^2) \|f'\|^2 \right\} +C\lambda^3\;;
   \end{eqnarray*}
since the first term at the \rhs can be made arbitrarily small, a
bound state exists as long as the curly bracket is negative. This
condition can be rewritten as
$$
\eta^2- 2\eta z +3z+K^2 \,<\,0\,, \qquad z:= \left(\pi\over a\right)^2
{\|f\|^2\over \|f'\|^2}\,,
$$
so it can be satisfied provided $\,z^2\!-3z-K^2>0\,$, which
requires in turn
$$
z\,>\, {3+\sqrt{9+4K^2}\over 2} \;;
$$
substituting for $\,K\,$ from (\ref{K explicit}), we arrive at
the condition (\ref{existence condition}).

It remains to find an upper bound to the eigenvalue
$\,\epsilon(\lambda)\,$. We put $\,\eta=z\,$, so
$$
(H_\lambda\psi,\psi) -\,\left(\pi\over a\right)^2 \|\psi\|^2
\,=\, \|\psi_x\|_{L^2(\Omega_b^c)} -\lambda^2 \|f\|^2 d_1
+\OO(\lambda^3)\,,
$$
where $\,d_1:= \left(\pi\over a\right)^2
(z\!-\!3\!-\!z^{-1}K^2) >0\,$. Outside $\,\Omega_b\,$ we choose
the trial function as follows:
$$
\psi(x,y)\,:=\; \left\{ \begin{array}{lll} e^{-\kappa|x\mp
b|}\chi_1(y) & \;\dots\; & \pm x>b \\ \\ 
0 & \;\dots\; & {\rm otherwise} \end{array} \right.
$$
where $\,\kappa=\,{1\over 2}\,\lambda^2 d_1\|f\|^2$.
In that case $\,\psi\,$ clearly belongs to the form domain of
$\,H_\lambda\,$ and 
$$
{(H_\lambda\psi,\psi)\over \|\psi\|^2} -\,\left(\pi\over a\right)^2
\,\le\: {-\,{1\over 2}\,d_1\|f\|^2\lambda^2 \over 2b-
{2\over \lambda^2 d_1\|f\|^2} -C\lambda}\,=\, 
-\,{1\over 4} d_1^2\|f\|^4 \lambda^4 +\OO(\lambda^5)\,,
$$
which yields the other bound of (\ref{asymptotics}).

\section{Concluding remarks}

The result discussed here is certainly not optimal. To illustrate
this, let us ask what is the maximum value of $\,{a\over b}\,$ for
which the sufficient condition (\ref{existence condition}) may be
satisfied.  Since the function $\,f\,$ satisfies the condition
(\ref{zero mean}) and $\,f(b)=0\,$, even the inequality
(\ref{embedding 2}) cannot be saturated for it; minimizing over the
smaller class, we obtain
$$
\inf\, {\|f'\|^2\over \|f\|^2}\,=\, \left(\pi\over b\right)^2\,.
$$
Hence examples of critical--tube shapes satisfying the condition
(\ref{existence condition}) exist if
$$
{a\over b}\,<\, \sqrt{2\over 3+\sqrt{9+4K^2}}\;\approx\, 0.506\;;
$$
on the other hand, the condition (\ref{nonexistence condition})
excludes the existence of bound states for the ratio $\,{a\over b}\,>
2.309\,$.  
\vspace{5mm}

\noindent
{\em Acknowledgments.} S.V. is grateful for the hospitality
extended to him in the Nuclear Physics Institute, AS, where this work
was done. The research has been partially supported by the Grants AS
No. 148409 and GosComVUZ 5--2.1--50.

\end{document}